\documentclass[a4paper,12pt]{article}
\usepackage{amsmath,amsthm}
\usepackage[margin=1.2in]{geometry}
\def\={\!=\!}

\usepackage[usenames]{xcolor} 

\def\={\!=\!}

\def\|{{\Vert}}

\newcommand{\noi}{\noindent}

\newcount\rem
\def\Remark{\skip\noi{{\bf{Remark \number\rem.}}} \advance\rem by 1}

\title{A measure theoretic  paradox from a continuous colouring rule\footnote{dedicated to Robert Aumann and Jan Mycielski}}
\author{Robert Simon, Grzegorz Tomkowicz}

\begin{document}
\begin{titlepage}
\maketitle
\thispagestyle{empty}

\begin{abstract} 
Given a probability space $(X, {\cal B}, m)$, measure preserving transformations $g_1, \dots , g_k$ of $X$, and a colour set $C$,  a colouring rule 
 is a way to colour the space with $C$ such that the colours  allowed for a
point $x$ are determined  by that point's location and the colours of the finitely  
 $g_1 (x), \dots , g_k(x)$  with $g_i(x) \not= x$ for all $i$ and  almost all $x$. 
 We represent a colouring rule as a correspondence $F$ defined on $X\times C^k$ with values in $C$. A function 
  $f: X\rightarrow C$ satisfies  the rule at $x$ if $f(x) \in F( x, f(g_1 x), \dots , f(g_k x))$.    
   A colouring rule is paradoxical if it can be satisfied
in some way almost everywhere with respect to $m$, but not in {\bf any} way that is measurable with
respect to a finitely additive measure  that extends the probability measure $m$ defined on ${\cal B}$ and for which the 
finitely many 
transformations $g_1, \dots , g_k$  remain measure preserving. Can a colouring rule be paradoxical if both $X$ and 
  the colour set $C$
 are convex and compact sets 
 and the colouring rule says if $c: X\rightarrow C$ 
  is the colouring function then  the  colour $c(x)$   must lie ($m$ a.e.)  in $F(x, c(g_1(x) ), \dots , c(g_k(x)))$ for 
 a non-empty upper-semi-continuous convex-valued correspondence  $F$ defined on 
  $X\times C^k$?   The answer is yes, and we present such an example. 
We show that this result is robust,  
 including  that any colouring that     
  approximates   the correspondence by  $\epsilon$ for small enough positive 
  $\epsilon$  also   cannot be measurable in the same finitely additive way.  Because non-empty 
  upper-semi-continuous convex-valued 
   correspondences on Euclidean space  can be approximated by continuous functions, there are paradoxical
     colouring rules that are defined by continuous functions.

\end{abstract}
\vskip.1cm

\noi Robert Simon \newline London School of
Economics, Department of Mathematics\newline Houghton Street,
London WC2A 2AE\newline e-mail: r.s.simon@lse.ac.uk\newline
\vskip.1cm
 
\noi Grzegorz Tomkowicz\\
Centrum Edukacji $G^2$,
ul.Moniuszki 9\\
41-902 Bytom,
Poland\\
e-mail: gtomko@vp.pl

\end{titlepage}



\section{Introduction}

In~\cite{ST1}, we introduced  colouring rules.    We demonstrated  several 
paradoxical colouring rules and proved that if there are finitely many colour
classes and the measure preserving transformations are invertible, 
  then any colouring of a paradoxical colouring rule has colour classes
that jointly, with the measure preserving transformations and the 
 Borel sets, define a {\em measurably $G$-paradoxical decomposition } (for the group $G$ generated by 
  the measure preserving transformations),  
  by which we mean the existence of two measurable sets of different
measures that are $G$-equidecomposable (see~\cite{ST1}, Thm. 1). 

In the conclusion of~\cite{ST1}, we asked whether a colouring rule could be
paradoxical if the colour classes belonged to a finite dimensional convex set
and the colouring rule was defined by an upper-semi-continuous convex-valued
non-empty correspondence, as described above. We call such colouring rules {\em probabilistic}
colouring rules. This means, among other things, that the choosing of colours
could be according to a maximisation or minimisation of a continuous and affine evaluation
of options, with indifference between two options implying indifference between
all of their convex combinations.           
        
        Our main inspiration is the question whether measure theoretic paradoxes, such as the Banach Tarski Paradox, 
         have any applications to areas  beyond mathematics, such as physics or economics.  A colouring rule could 
          represent natural forces, and the lack of any measurable solution could represent 
           a radical  inability to predict their  behaviour.           
A particular  inspiration is the widely held belief in economic theory that
although one cannot always accomplish optimisation goals through behaviour that
is measurable with respect to a countably additive measure,  one can do so  with some finitely additive option.  The
problem with this belief  is that there may be  knowledge structures to
the optimisation that cannot be altered when extending to a finitely additive
measure. If those knowledge structures are defined through the use of ergodic
operators, measure invariance of those ergodic operators may be required.

             The Brouwer Fixed Point Theorem is relevant to probabilistic colouring rules. 
          If   $X$ is compact and the correspondence is independent of the location of $x$, meaning that 
        the correspondence    $F$ is defined entirely on $g_1(x), \dots , g_k(x)$, then the fixed point theorem shows there exists 
            a constant colouring function satisfying the colouring rule, hence it cannot be paradoxical. This is done by 
             mapping $C$ to $k$ copies of $C$ through $k$ copies of the identity,  and then following the colouring 
              rule back down to $C$ via the correspondence $F$.
               Therefore we have to consider colouring rules that are dependent on the 
               location in the space. 
       
       In Simon and Tomkowicz [ST] we demonstrated a probabilistic colouring rule with a one dimensional continuum 
        of colours such that after the correspondence is  approximated by any $\epsilon$ for small enough positive $\epsilon$ 
          the colouring rule  still had  no Borel measurable solution.  
          This  colouring rule was  paradoxical only in the sigma-additive, not in the  finitely additive, sense.

In the next section we describe the probabilistic colouring rule and show that it is paradoxical. In the
third section  we look at approximating the colouring rule  and apply it 
  to economics.   In conclusion we consider
related problems.
 
\section{A Probabilistic Paradoxical  Colouring Rule} 
 
 Let $T_1$ and  $T_2$ be  two non-invertible generators (each  generating a semi-group
  isomorphic to ${\bf N}$). We assume that $T_1^3 T_2 = T_2 T_1^3$ and there are no other relations.  Let $G$ be the   
  semi-group  generated  by $T_1$ and $T_2$.  
   Let $X$ be the set $\{ 0,1\}^G$. We extend $X$ to $X' = X \times \{ a,b,c\}$ and let the  
    symmetric group $S_3$ act on the three elements $\{ a,b,c\}$. We assume that 
    $S_3$ commutes with $G$ and define $G'$ to be the semi group so generated. We need the addition 
     of $S_3$ to define the colouring rule according to measure preserving transformations in $G'$. However it is 
      mostly the colouring of 
      $X\times \{a \}$ that matters.  
 
   For  any $x\in X$ and $g\in G$, $x^g$ stands for the $g$
coordinate in $x$. With $e$ the identity in $G$, the $e$ coordinate of $x$ is
$x^e$.  There is a canonical right semi-group action on $X$, namely
$g(x) ^h= x^{gh}$ for every $g,h\in G$. We use the canonical product topology on $X$. 
  For every cylinder determined by
particular  choices of $\{ 0,1\} $   
 we assign the probability $ (\frac {1} {2} ) ^k $ where $k$ is the number of those choices 
  determining the cylinder.   With this Borel
probability measure the semi-group $G$ act measure preserving on $X$. Any semi-group element acts measure preserving  on any cylinder, due to the cancellation law, and this can be extended to any Borel set through approximation via cylinders. 
Likewise we give $X'$ the Borel 
 probability measure where each element in $\{ a,b,c\}$ is given equal probability when paired with a 
  Borel set in $X$. The probability of 
  $A \times \{ a,b,c\}$ in $X'$ is given the same probability at that given to  $A$ in  $X$, 
    In this way  $G'$ acts measure
  preserving on $X'$. With  $m$ the  canonical  Borel measure of $X$, let $m'$ be its extension to  $X'$.

   The subset of $X$ where $g\not= h$ implies that  $gx\not= hx$ is of Borel measure one. This follows 
    from the fact that the semi-group is countable. 
   Without loss of generality, we will be interested only in this subset, and we ignore the set of measure $0$ where 
    this doesn't hold.

The set of colours $C$  is $\Delta (\{1,2,3\}):= \{ p \ | \ \forall i \ p_i\geq 0, p_1+ p_2 + p_3 =1\} $ where $\delta_i$ is perceived 
 to be all weight to the colour $c_i$.  We represent the extremal 
  colours $c_i$  modulo $3$, with $i=1,2,3$ rather than $0,1,2$.

Now we show how to colour each $x\in X'$ with a point in 
 the  simplex $\Delta (\{ 1,2,3\})$, with  $p=(p_1, p_2, p_3)  \in \Delta (\{ 1,2,3\})$ standing 
 for the weights given to the three colours. We colour according to the optimisation of a continuous function, and 
  the resulting optimal solutions are represented by a correspondence $F$, defined on $X\times C^2$, where each $x\in X$ has 
    two descendants $T_1 x$ and $T_2 x$.  The continuous functions are defined by matrix multiplications. 
   
    For every   $x\in X$ we  define   
  a matrix $A_x=  
  \begin{pmatrix} 1 & r^x_{1,2}  & r^x_{1,3}  \cr 
  r^x_{2,1}  & 1 & r^x_{2,3}   \cr r^x_{3,1} & r^x_{3,2} & 1 \end{pmatrix} $
     such that  $|r^x_{i,j}| \leq \frac 1 {100}$ for 
   every choice of $x$ and $i,j$. 
   The $r^x_{i,j}$ are continuous functions of $ x$ in $X$, and  are chosen such that 
    for every  $p\in \Delta (\{ 1,2,3\}) $ the 
     set of $x$ such that only one row   maximises $A_x p$ is a set of Borel  measure one. The following  is one way to do that.   There are six entries of 
      the $r^x_{i,j}$ to determine. 
      We  place the elements of $G$ into six infinite ordered collections. For each $l=1,2,\dots$ and $g_l$ the $l$th group elements corresponding $r^x_{i,j}$ we  
       define $r^x_{i,j}:= \frac 1 {100}  \sum_{l=1} ^{\infty}  2^{-l}  x^{g_l}$. where the structure of $X$ requires 
        that $x^g\in \{ 0,1\}$ for all $g\in G$.

  We colour a point $(x,b)$ by any convex combination of the $c_i$ with the  $i=1,2,3$ whose rows maximise  $A_{T_2 T_1 T_2 (x)} p$, 
   where $p$ is the colour given to $(T_1 x , a)$. 
   
   We colour a point $(x,c)$ by any convex combination of the $c_i$ with the  $i=1,2,3$ whose rows maximise  $A_{T_1 T_2 T_1 (x)} p$, 
   where $p$ is the colour given to $(T_2 x , a)$.

  The complex part of the rule is how to colour a point $(x,a)$. 
    
    We define two  three  dimensional   matrices $B_0$ and $B_1$. The matrix 
         $B_0$ is used when $x^e=0$ and the matrix $B_1$ is used when $x^e=1$. Both matrices 
          have  entries $(b_{i,j,k}) $  only in  $\{ 0,1\}$ such that for each pair of columns 
     $j,k$ there is only one non-zero $b_{i,j,k}$ and it is equal to $1$. 
 There are two sets 
     of columns and one row; the columns $j$ and $k$ correspond to the weights given to the colours $c_j$ and $c_k$ of   
      $(x,b)$ and  $(x, c)$ respectively.

         A row $i=1,2,3$ of $B_1$ or $B_0$  is evaluated in the following way. If $p\in \Delta (\{ 1,2,3\})$ is the colour 
         given to $(x,b)$ and $q\in \Delta (\{ 1,2,3\})$ is the colour 
         given to $(x,c)$ then the $i$th row is given the value $\sum_{j=1}^3 \sum _{k=1}^3 p_j q_k b_{i,j,k}$. The colouring 
          rule requires that any convex combination of rows is chosen that maximise this row evaluation.  

        To define these matrices, we determine when the $b_{i,j,k}$ entry is $1$.    
        The $j$th column represents the an extremal colour of $(x,b)$ and the $k$th column an extremal  colour of $(x,c)$.  
        
        The matrix $B_1$, the case of $x^e=1$, is easy to define. 
         The $b_{i,j,k}$ entry is $1$ if and only if $i=j+1$.
        
        The matrix $B_0$, the case of $x^e=0$, is more complex. 
        
        If  $k\not= 1$, and $j\not=3$, then the $b_{i,j,k}$ entry is $1$ if and only if $i=j$.
             
             If  $k\not= 1$,  then the $b_{i,3,k}$ entry is $1$ if and only if $i=1$.   
             
                    If  $j\not=3$, then the $b_{i,j,1}$ entry is $1$ if and only if $i=j+1$. 
                  
                    The $b_{i,3,1}$ entry is $1$ if and only if $i=3$. 
                    
                   This completes the definition of the correspondence $F$. We will also refer to $F$ as the colouring rule.   Notice that when $x^e=0$ then 
                    the matrix $B_0$ has similarity to the conditions of the Hausdorff paradox, requiring that if
                     the colour of $x$ is an advancement by one  on the colour of $T_1 x$ and $x^e=0$ then
                      the two  points $x$ and $T_2 x$ are coloured differently, one of these two points is coloured 
                        $\{ c_1\}$  and the other takes a colour in    $\{ c_2, c_3\}$.

    When $(x,a)$, $(x,b)$ or $(x,c)$  gives all weight to an 
    extremal colour, namely $\delta_{c_i}$ for some $i\in \{ 1, 2, 3\}$, then it is called 
  {\em pure}.

 A colouring $c: X'\rightarrow C$ {\em satisfies} the colouring rule $F$ if the rule holds almost everywhere with respect to the 
  probability distribution $m'$.   
    \vskip.2cm 
  
   {\bf Lemma 1:} A colouring that satisfies the colouring rule $F$ is pure almost everywhere with respect to the Borel measure 
   $m'$. 
   
   {\bf Proof:} There is no relation between $T_2 T_1 T_2$ and $T_1$, likewise between 
    $T_1 T_2 T_1$ and $T_2$.  Therefore for every choice of $x$ the 
     matrices at  $(T_1 ^{-1} x,b)$ at  $(T_2^{-1} x,c)$  defining the colouring rule $F$  are the matrices  $A_y$ for all $y\in X$.
      The conditional probability on those matrices 
       is the same as the distribution on $X$. It follows by the construction of these matrices 
       that regardless of the colour  $p$ at $(x,a)$  the set of colours at $(T_1^{-1} x, b) $ and $(T_2 ^{-1}x, c)$  that 
        are not pure with respect to $m'$  is a  subset of conditional  measure $0$. By 
           the definition of conditional probability, the conclusion follows for all $(x,b)$ and $(x,c)$.  And therefore 
            it follows for $(x,a)$ also.  
           \hfill q.e.d.

              \vskip.2cm 
       Lemma 1 allows us to perceive a colouring of $X'$ as being primarily a colouring 
        of $X \times \{ a\}$. With purity, the colours of $(x,b)$ and $(x,c)$ are merely conveying to 
         $(x,a)$ the colours of $(T_1 x, a)$ and $(T_2 x, b)$ in a way that those two colours are 
          determining the colour of $(x,a)$. 
           From now on, by a colouring of $X$ we mean a colouring of $X\times \{ a\}$, where 
         by the colour for $x$ we mean the colour for $(x,a)$.  If a point is coloured purely 
          with   $\delta _i$ we will also write that it is coloured with $c_i$.   Our main aim is to prove the following theorem:

        {\bf Theorem 1:}   For any 
         finitely additive $G$-invariant  measure $\mu$  on $X'$ extending $m'$  there exists no 
          colouring $c:X'\rightarrow C$ that is  $\mu$ measurable and satisfies the colouring rule $F$.
          \vskip.2cm 
          
       {\bf Definitions:} Semi-group elements 
         $g_1, g_2, \dots , g_k\in G$ are called {\em independent} if there are no relations between them.  Two 
          points $x, x'$ are called {\em twins} if $T_1 x = T_1 x'$ and $T_2 x = T_2 x'$, meaning 
           that the differ only by $x^e \not= x^{x'}$. A point $x\in X$ is 
          coloured {\em randomly} if  $x$, the twin of $x$,  $T_1 x$, and  $T_2 x$ are all  pure  and $x$ and its twin are coloured 
           differently (meaning that for the twins $x, x'$  the colour of $T_1x$ is advanced by one to define the colour of $x$ where 
            $x^e=1$  and the colour of $T_1x$ is not advanced  to define the colour of $x'$ where $(x')^e =0$).
             A  colouring  that is measurable with respect to any  $G$-invariant finitely additive measure 
            $\mu$ is called {\em measurable}. 
          \vskip.2cm 
        
        {\bf Lemma 2:} If there is a positive measure  
         of points that are coloured randomly for a measurable colouring satisfying the colouring rule $F$, and $g_1, g_2, \dots, g_k$ are independent, 
          then the probability of $\{ x \ | \ g_i x \mbox { is coloured }  c_{n_i}, i=1, 2, \dots, k\} $  is the product 
           $\prod_{i=1} ^k q_{n_i} $ where $q_n $ is the probability of $\{ x \ | \ x \mbox { is coloured }  c_n\} $ for 
            $n=1,2,3$. 
              \vskip.2cm

         {\bf Proof:} First we show that if there is a positive probability of 
        random colouring, then from the stochastic matrix representing 
        the transition of the distribution of   colours of $x$ to the colouring of   $T_1^{-1}x$ 
        there is a unique eigenvector in $\Delta ( \{ 1,2,3\})$  with the eigenvalue  
        $1$ and the other eigenspaces correspond to eigenvalues  with norms less than $1$. Since by measure invariance we can assume that 
         the probability distribution of the colours is the same at $x$ and $T_1^{-1}x $ for the collection of all $x\in X$, this 
          implies that the distribution is determined by this eigenvector. 
        
         Let $p,q,r$ be  half the probabilities for the  random colouring  conditioned on the colours 
          $c_1, c_2, c_3$ respectively, with the half referring to the  non-advancement 
          of the colours  which happens when $x^e=0$.  We assume that at least 
         one of the $p,q,r$ are positive and none are greater than $\frac 12$.  
         The stochastic matrix in question is
         $$\begin{pmatrix} p & 0 & 1-r \cr 1-p & q &0 \cr 
          0 & 1-q & r\end{pmatrix}.$$ 
          It has the characteristic polynomial $(p-x) (q-x) (r-x) + (1-r) (1-q) (1-p) = 
          (1-r) (1-q) (1-p)+ pqr - (qr+ qp+ pr) x +  (p+q+r) x^2- x^3$. The degree one polynomial   $1-x$ divides 
           this characteristic polynomial, leaving  the second degree polynomial 
            $1-r-p-q + rp+ rq+ qp + (1-r-p-q)  x +   x^2$ as a factor. Letting $b=1-r-p-q$ and $d=  rp+rq+qp$ 
             we have the roots $\frac { -b + \sqrt {b^2 - 4 (b+d)  }} 2$ and  $\frac { -b - \sqrt {b^2 - 4 (b+d) } } 2$.
            Assuming that $b^2-4(b+d)$ is not positive,  
            the norm squared of these roots is equal to 
             $b+d$, which is less than $1$ because one of $p,q,r$ must be positive 
              and if, for example,   $r>0$ then $r> rp$.  
              $b^2- 4b$ cannot be positive if $b$ is not negative.  
              Since $b\geq -\frac 12$ and $\frac 12 + \sqrt {2+\frac 14} =2$, 
               the only way for the norm of $\frac { -b - \sqrt {b^2 - 4 b  }} 2$ or $\frac { -b + \sqrt {b^2 - 4 b  }} 2$ to reach
                $1$ is if $b=-\frac 12$ and therefore $r=p=q = \frac 12$.  But then   
                 the roots are really $\frac { -\frac 12  + \sqrt {-\frac 34  }} 2$ and  $\frac { -\frac 12 - \sqrt {-\frac 34 } } 2$, 
                  with norms of $\frac 12$. 
              
             Due to the lack of any relation between the $g_i$, with respect to the Borel probability distribution
               whether the  $(g_i x)^e$ are equal to $0$ or $1$ are independent 
              choices over all the $x\in X$, hence also with 
               any finitely additive measure extending the Borel measure.  Assuming invariance of joint distributions of colour combinations before and 
               after  applying $T_1^{-1}$, using induction on the number of the $g_i$, 
                 and that the independence of the $g_i$ implies the independence of the $T^n_1g_i$ for all positive $n$, 
                   we complete the proof with the following claim: 
              
              {\bf Claim:} Assume that there are two one-stage 
              stochastic processes on two finite sets $S$ and $T$ respectively and a stochastic process 
               defined on $S\times T$  such that the transitions are defined independently by  transitions on $S$ and $T$. 
                Assume for each of the $S$ and $T$ processes   that there is only one unique  invariant distribution/eigenvector 
                 corresponding to an eigenvalue 
                  of norm 
                 $1$. Furthermore assume 
                  that for each $s\in S$ and $t\in T$ that there is a positive probability that there is a transition to something other than 
                   $s$ and $t$ respectively. Then there is one  unique invariant joint distribution on $S\times T$ defined by 
                    the independent distributions on $S$ and $T$ respectively. 
                    
                    {\bf Proof of Claim:} Let $(s,t)$ be any pair of states in $S\times T$; we want to prove that any invariant probability 
                     for $(s,t)$ is the same as $ab$ where $a$ is the invariant probability for $s$ and $b$ is the invariant probability for $t$. 
                      Let $q$ be an invariant probability for $(s,t)$; we want to show that $q=ab$. Let $l_a$ be the probability, conditioned 
                       on the state being at $a$, of  leaving  $a$ on the next stage. Let $l_b$ be the same probability for the state $b$. 
                        Let $r_a$ be the probability of not being at $a$ and returning to $a$ on the next stage. Let $r_b$ be the 
                         same probability for the state $b$. We can calculate $q$ by $q= q- l_a q - l_b q + l_a l_b q + (1-a) r_a (1-l_b) b+ 
                          (1-b) r_b (1-l_a) a + (1-b) r_b (1-a) r_a $. If these distributions on $S$ and $T$ are  
                           invariant independently, then $\overline q=ab$ is another solution for an invariant probability for $(s,t)$.  Therefore 
                            we  can also write   $ab= ab- l_a ab - l_b ab + l_a l_b ab + (1-a) r_a (1-l_b) b+ 
                          (1-b) r_b (1-l_a) a + (1-b) r_b (1-a) r_a $, the same formula but with   $q$ replaced by $ab$. But then we can 
                           write $0 = q (l_a + l_b - l_a l_b) = ab (l_a + l_a - l_a l_b)$. With $l_a>0$ and $ l_b>0$
                             we have $l_a + l_a - l_a l_b>0$ and $q=ab$. \hfill q.e.d.
              \vskip.2cm   
              
         Notice that if $h= jg$ for semi-group elements $g,h, j$ then the choice of $gx$ for any $x\in X$  will determine the $hx$, including 
          the value of $(hx)^e$, and therefore  
           the $(gx)^e$ and $(hx)^e$ are dependent for the various choices of $x\in X$ in the maximal way that two variables 
            can be dependent. 
           If  $kh=jg$ for some semi-group elements $g,h,j,k$, then starting at 
           $kh(x) = jg(x) $  the dependence from  $h$ and $ jg$  can have residual influences on the  joint  colour distribution  at the set  
           $\{ (gx, hx) \ | \ x\in X\} $ before they reach their limit distributions.   This residual dependence implies
            that   the   joint colour distribution at  $\{ (T_1x, T_2x) \ | \ x\in X\} $  should have some dependence 
             coming from $T_1^3T_2 = T_2 T_1^3$. However  from a different 
             perspective they  should be independent!  This  contradiction 
             drives our proof of Theorem 1.  
            
           \vskip.2cm   
       
       {\bf Lemma 3:} There is no measurable colouring satisfying the colouring rule $F$ such that there is a positive 
        measure of points that are coloured randomly. 
         \vskip.2cm 
         
         {\bf Proof:} 
         Let $C_z $ be the stochastic matrix 
                $C_z = \begin {pmatrix} \frac {1-z} 2& 0 & \frac {2-z} 2 \cr 
                 \frac {z+1} 2 & \frac {1 -z} 2 & 0\cr 
                 0 & \frac {z+1} 2&  \frac z 2 \end {pmatrix}$.  Given an independent distribution of colours for the pair 
                  $T_1 x$ and $T_2 x$, the matrix 
                  determines the distribution of colours at $x$ where $z$ is the probability for the colour $c_1$ at $T_2x$.

         First consider the pair of independent elements $T_2 T_1$ and $T^2_2$ and let $(y_1, y_2, y_3)$ be the global 
          probabilities for the colours $c_1, c_2, c_3$, respectively. From Lemma 2 the joint  distribution of 
           the colours of $T_2 T_1 x$ and $T_2^2 x$ are that determined by the products of the  $y_i$, meaning 
            that the probability of  $T_2 T_1 x$ and $T_2^2 x$ coloured $c_i$ and $ c_j$ respectively are 
             $y_i \cdot y_j$. But as the colouring is measurable, using that $T_2$ is measure preserving, 
              we have to assume that the probability 
              of the $T_1 x$ and $T_2 x$ coloured $c_i$ and $c_j$ respectively is also 
               $y_i \cdot y_j$.  Applying the matrix $C_z$ with $z= y_1$  for $C_{y_1}  (y_1, y_2 , y_3)^t = 
                (y_1, y_2, y_3)^t$ we get the following three equations with three variables: 
         
         $$y_1 = \frac 12 (y_3) + \frac 12 ((1-y_1) y_1 + (1-y_1) y_3) $$
         $$y_2= \frac 12 (y_1) + \frac 12 (( 1-y_1) y_2 + y^2_1)$$
         $$y_3 = \frac 12 (y_2) + \frac 12 (y_1 y_3 + y_1 y_2).$$
         
         This solves to $y_1= 3-\sqrt 7$, $y_2= 3-\sqrt 7$, and $y_3 = 2 \sqrt 7 -5$. These are the global probabilities 
          for the colours $c_1, c_2, c_3$ respectively. 
         Approximately this is the triple $( .35425, .35425, .2915)$. The reason for a smaller probability for the colour $c_3$ 
          results from the tendency to move away from this colour with  greater probability than toward it. 
          
          Second  consider the seven semi-group elements. 
 $$g_1:= T_2^4 T_1 T_2^3, $$
$$g_2:= T_2^4 T_1T_2^2T_1T_2,$$
$$g_3:=T_2^4 T_1T_2^2T_1^2T_2, $$
$$g_4: =T_2^4 T_1T_2^2T_1^3, $$
$$g_5 := T_2^4 T_1T_2T_1^4, $$
$$g_6 := T_2^4 T_1T_2^2T_1^2, $$
$$g_7:= T_2^4 T_1T_2^2T_1.$$ We show that they are independent, meaning that there is no relations between 
 them.  Make the revision 
  that $T_1$ is invertible and $T_1^3$ is the identity. If the resulting elements still  have no relations between them, 
  there was no relation  before this revision. We can now proceed with the   assumption that there are no relations other than 
   $T_1^3= e$. 
  
  We have the reductions $$h_1:= T_2^4 T_1 T_2^3, $$
$$h_2:= T_2^4 T_1T_2^2T_1T_2,$$
$$h_3:=T_2^4 T_1T_2^2T_1^2T_2, $$
$$h_4: =T_2^4 T_1T_2^2, $$
$$h_5 := T_2^4 T_1T_2T_1, $$
$$h_6 := T_2^4 T_1T_2^2T_1^2, $$
$$h_7:= T_2^4 T_1T_2^2T_1.$$ 
 Notice that any combinations of the $h_i$ are  separated by  $T_2^4$ or a higher power of $T_2$, so that we can identify when one element ends and the next begins. Removing these bookends of $T_2^4$, $T^5_2$,  $T^6_2$, or $T^7_2$, we are left with  
 the interiors $ T_1, T_1T_2^2T_1, T_1T_2^2T_1^2, 
   T_1T_2T_1.$ All but the $T_1 T_2 T_1$ (from $h_5$) are repeated twice. The two uses of 
   $T_1$ are distinguished by the powers of $T_2$ succeeding them ($T_2^7$ vs $T_2^6$). 
   Likewise the two uses of $T_1T_2^2T_1$ can be distinguished 
  by the power of $T_2$ succeeding them ($T_2^5$ vs $T_2^4$) and the same is true 
   of $T_1 T_2^2 T_1^2$ ( $T_2^5$ vs $T_2^4$). It follows that 
    all seven elements $g_1, \dots, g_7$  also must  be independent. 
 
 By Lemma 2 the probability of 
 $\{ x \ | \ g_{i}  x\mbox { coloured } c_{n_i} \} $ for any choice of $n_1, \dots , n_7$ is equal to 
  the product $\prod _{i=1} ^7 y_{n_i} $. Due to the measure preserving property of $T_2^4 T_1 T_2$, 
   we can assume that the probability of 
 $\{ x \ | \ g'_{i}  x\mbox { coloured } c_{n_i} \} $ for any choice of $n_1, \dots , n_7$ is equal to 
  the product $\prod _{i=1} ^7 y_{n_i} $ for the elements   $$g'_1:= T_2^2,$$
$$g'_2:= T_2T_1T_2,$$
$$g'_3:=T_2T_1^2T_2, $$
$$g'_4: =T_2T_1^3, $$
$$g'_5 := T_1^4, $$
$$g'_6 := T_2T_1^2, $$
$$g'_7:= T_2T_1.$$ 
          
          Now consider the defining condition that $g'_4 x = T_1^3 T_2 x = T_2 T_1^3 x $ for all $x$. 
          For each of the three   colour choices  at $T_1^3 T_2 x$ we determine  
           the distribution on the colours of $x$ by building up to a joint distribution for $T_1 x$, $T_2 x$, and $x$ through 
            the two pathways, the left pathway going through three applications of $T_1^{-1}$ and the right 
            pathway going first through $T_2^{-1}$ followed by two applications of $T_1^{-1}$. Due to 
            the above mutual independence of colour distributions of the $g'_ix$ and their locations 
             in relation to the  two pathways we can use the matrices $C_z$ to 
             determine these distributions, including the last step of determining the colour 
              distribution of $x$, conditioned on the colour of $T_1^3 T_2 x$.  
                  There are three matrices most relevant to our calculations, 
                  $C_1$, $C_0$, and $C_{y_1}$. 
                  
                  The distribution of colours for $T_2 x$, conditioned on a choice of 
                   colour for $T_1^3 T_2$ is determined by 
                   $$C_{y_1} ^3 = \begin{pmatrix} \frac 1 8 (50\sqrt 7 - 129) & \frac 14 (23-8\sqrt 7) & \frac 1 8 ( 19\sqrt 7 - 49) \cr 
                    \frac 18 (216- 81\sqrt 7 ) & \frac 18 (50 \sqrt 7-129) & \frac 14 (23-8\sqrt 7)\cr 
                     \frac 18 (31\sqrt 7 -79) & \frac 18 (91-34 \sqrt 7) & \frac 18 (11-3\sqrt 7) \end{pmatrix}.$$ To determine the distribution 
                     of colours for $T_1 x$ we need to calculate $ C_{y_1} ^2  C_0 (y_1, y_2, y_3)$ and 
                     $ C_{y_1} ^2  C_1 (y_1, y_2, y_3)$, the former for  what happens when  $T_1^3 T_2 x$ is  coloured $c_2$ or $c_3$ 
                      and the latter for $T_1^3 T_2 x$ coloured $c_1$. For the former we get 
                        $$\begin{pmatrix} \frac 14 (11-4\sqrt 7) & \frac 14 (5\sqrt 7 -11) & \frac 14 (\sqrt 7-1) \cr 
                         \frac 12 (6\sqrt 7-15) & \frac 14 (11- 4\sqrt 7) & \frac 14 (5\sqrt 7 -11) \cr 
                          \frac 14 (23-8\sqrt 7) & \frac 14 (4-\sqrt 7 ) & \frac 12 (8-3\sqrt 7) \end{pmatrix} \begin{pmatrix} 
                          \frac 32 \sqrt 7-\frac 72\cr  
                          3-\sqrt 7\cr  \frac 32 - \frac {\sqrt 7} 2)\end{pmatrix}  =$$
                          $$\big(\frac 18 (117\sqrt 7 - 307) , \frac 14 (258-97\sqrt 7), \frac 1 8(77\sqrt 7-201\big).$$
                   For the latter we get $$ \begin{pmatrix} \frac 14 (11-4\sqrt 7) & \frac 14 (5\sqrt 7 -11) & \frac 14 (\sqrt 7-1) \cr 
                         \frac 12 (6\sqrt 7-15) & \frac 14 (11- 4\sqrt 7) & \frac 14 (5\sqrt 7 -11) \cr 
                          \frac 14 (23-8\sqrt 7) & \frac 14 (4-\sqrt 7 ) & \frac 12 (8-3\sqrt 7) \end{pmatrix} \begin{pmatrix} 
                           \sqrt 7 - \frac 5 2\cr 3-\sqrt 7 \cr   \frac 12   \end{pmatrix}  = $$ 
                            $$\big(\frac 18 (95\sqrt 7 - 248), \frac 18 (429-161\sqrt 7), \frac 18(66\sqrt 7-173) \big).$$ 
                           
                        For each of the three possibilities, $T_1^3 T_2 x = T_2 T_1^3 x$ coloured $c_1$, $c_2$, or $c_3$, we
                             do the calculations for the three colours at $T_1^3 T_2 x$ separately and recombine the results 
                             according to $(y_1, y_2, y_3)$. 
                             
                             For $T_1^3 T_2 x$ coloured $c_1$, we have to calculate $C_{\frac 18(50\sqrt 7 -129)}  
                               (\frac 18 (95\sqrt 7 - 248), \frac 18 (429-161\sqrt 7), \frac 18(66\sqrt 7-173) ^t =$
                               $$\begin{pmatrix} \frac 1 {16} (137-50\sqrt 7)& 0 & \frac 1 {16} (145-50\sqrt 7) \cr 
                                \frac 1 {16} (50 \sqrt 7 -121) & \frac 1 {16} (137-50\sqrt 7) & 0 \cr 
                                 0 & \frac 1 {16} (50 \sqrt 7 -121 ) & \frac 1 {16} (50\sqrt 7 -129)\end{pmatrix} 
                                  \begin{pmatrix} \frac 18 (95\sqrt 7 - 248)\cr  \frac 18 (429-161\sqrt 7)\cr  \frac 18(66\sqrt 7-173)\end{pmatrix} 
                                  $$ =
                                  $$\big(\frac 1 {128} (43635\sqrt 7-115411) , 
                                   \frac 1 {128} (178381-67402\sqrt 7) , \frac 1 {128} (23767\sqrt 7 -62842)\big)$$
                                
                               For $T_1^3 T_2 x$ coloured $c_2$, we have to calculate 
                                $C_{\frac 14 (23-8\sqrt 7) }  (\frac 18 (117\sqrt 7 - 307) , \frac 14 (258-97\sqrt 7), \frac 1 8(77\sqrt 7-201)^t=$ 
                              $$\begin{pmatrix} \frac 1 8 (8\sqrt 7 -19) & 0 & \frac 18 (8\sqrt 7 -15) \cr 
                               \frac 1 8 (27-8 \sqrt 7) &   \frac 1 8 (8\sqrt 7 -19) &  0\cr 
                                0 &  \frac 1 8 (27-8 \sqrt 7) & \frac 1 8 (23-8 \sqrt 7)  \end{pmatrix} \begin{pmatrix} 
                               \frac 18 (117\sqrt 7 - 307) \cr \frac 14 (258-97\sqrt 7)\cr \frac 1 8(77\sqrt 7-201) \end{pmatrix} $$ = 
                               $$\big( \frac 1 {32} (9856-3721\sqrt 7) ,\frac 1{64} ( 13429\sqrt 7- 35509)  ,\frac 1 {64} (15861- 5987\sqrt 7)   \big).$$
                               
                                 For $T_1^3 T_2 x$ coloured $c_3$, we have to calculate 
                                  $C_{\frac 18 (19\sqrt 7-49) }(\frac 18 (117\sqrt 7 - 307) , \frac 14 (258-97\sqrt 7), \frac 1 8(77\sqrt 7-201)^t =$
                                  $$\begin{pmatrix} \frac 1{16} (57- 19\sqrt 7) & 0 & \frac 1 {16} (65-19 \sqrt 7) \cr 
                                   \frac 1 {16} (19 \sqrt 7-41) &   \frac 1{16} (57- 19\sqrt 7) & 0\cr 
                                    0 &  \frac 1 {16} (19 \sqrt 7-41) &  \frac 1 {16} (19 \sqrt 7-49) \end{pmatrix} \begin{pmatrix} 
                               \frac 18 (117\sqrt 7 - 307) \cr \frac 14 (258-97\sqrt 7)\cr \frac 1 8(77\sqrt 7-201) \end{pmatrix} =$$
                               $$\big( \frac 1 {64} (10663\sqrt 7- 28183), \frac 1 {64}
                                (41681-15746 \sqrt 7), \frac 1{64} (5083\sqrt 7 - 13434)\big).$$
          
          Combining these three vectors according to the distribution $(3-\sqrt 7, 3-\sqrt 7, 2\sqrt 7-5)$ of the three colours 
           at $T_1^3 T_2 x$ we get a different distribution for the colours at $x$, namely 
            $(\frac 1{64} (75588- 28561\sqrt 7), \frac 1 {128} (95189\sqrt 7 -251801) , \frac 1 {128} 
            (100753- 38067\sqrt 7))$, a contradiction. In decimals this corresponds approximately to\newline  
             $(.35464, .35486, .29050)$, (different from  the original $( .35425, .35425, .2915)
              =(3-\sqrt 7, 3-\sqrt 7, 2\sqrt 7 -5 )$).

                        \hfill q.e.d.  \vskip.2cm 
                
                The source of the discrepancy in the final distribution at the top 
                 comes from  different starting distributions at $T_2 T_1^3 x$ 
                 and $T_1^3 T_2 x$. If 
                 there was no commuting of $T_1^3$ and $T_2$, there would be two different points and the probability 
                  of a  start  of  
                  $c_i$ at $T_1^3 T_2x$ and $c_j$ at $T_2 T_1^3 x$  
                 would be  $y_i y_j$   with $(y_1, y_2, y_3) = (3-\sqrt 7, 3-\sqrt 7, 2\sqrt 7 -5)$. But with 
                  the commuting there are  no mixed starts.    
                
                Finally with the above lemmas, we can prove Theorem 1. 
                             \vskip.2cm 
                                      
                                    {\bf Proof of Theorem 1:} By the three lemmata,  a measurable colouring satisfying  the colouring 
                                     rule must have  no positive subset of points coloured randomly. 
                                     That means, with any probability approaching $1$,  that the colours cycle through with $x$ coloured 
                                       $c_{i+1}$ whenever $T_1$ is coloured $c_i$.  
                                      So with the measure preserving 
                                       application of $T_1$ we see that $\frac 13$ of the space is coloured $c_1$. 
                                       On the other hand, by $T_2$ being measure preserving, 
                                as with the Hausdorff paradox, 
                                this  implies that the probability of the space coloured $c_1$ must be  
                                 arbitrarily close to $\frac 23$, a contradiction.  \hfill q.e.d. 
              \vskip.2cm 
              
          There is of course a non-measurable colouring of the space. Start at some  $x$. 
           Let the colours cycle through $c_i$ by the repeated 
           application of $T_1^{-1}$, increasing the colour by one, and $T_1$, decreasing the colour by one. Choose some   $y$ already  coloured this way and  
              colour $T_2 y$ so that 
           that $T_2 y$ is coloured $c_1$ if and only if 
            $ y$ is not coloured $c_1$. Extend this choice for $T_2 y$ with  the repeated application of $T_1^{-1}$ and $T_1$  
             to $T_2 y$. Notice that   
             the commuting of $T_1^3$ and $T_2$ does not get in the way of this pattern. This can be continued, but there   
              is  a general solution to the existence of a non-measurable colouring which 
             includes this type of  colouring. This general solution is 
              presented in the next theorem. It demonstrates the other half of the argument that $F$
                is a paradoxical colouring rule, that 
               the colouring rule $F$ can be satisfied.

   Let $X^*$ be the $x$ in  
     $X$ where $g\not= h$ implies that  $gx\not= hx$. 
   \vskip.2cm
       
           {\bf Theorem 2:} There are pure  colourings of $X^*$ that satisfy the colouring rule $F$. 
             \vskip.2cm 
             
          {\bf Proof:}  Let a subset $A$ of $X$ be called {\em  closed} 
          if whenever  $T (x) \in A$ and $T_2 (x) \in A$ then also $x$ is in $A$.   Let 
            $\overline A$ stand for the closure of $A$. 
            A subset $A$ of $X$ is called {\em pyramidic} if whenever $x\in A$ then the $G$ orbit $Gx$ is in $A$.    Notice that the 
            closure of a pyramidic set is pyramidic.         
             
             For any pyramidic set $B$ and any pure  
              colouring of $B$ consistent with the colouring rule 
                there is a deterministic way to colour the closure  $\overline B$ according to the colouring rule.  
         We define a partially ordered set on the pairs $( B, c)$ where $B$ is a pyramidic and closed set 
          and $c$ is a colouring of $ B$  according 
          to the rule.   We say that  $(B,c) \geq (B' ,c')$ if $B$ contains $B'$ and $c$ restricted to $B'$ is $c'$.

          With Zorn's lemma there is a maximal element to any tower of the partial order. We  show 
          that  it is not possible for $(B,c)$ to be maximal however there is some $x\in X^*$ such that $x\not\in B$. 
          
          Consider the orbit $Gx$ and the sequence $A_l$ where 
           $A_l$ is the subset of $Gx$ such that $y\in A_l$ if and only if $y = gx$ with $g$ a word of length 
             $l$ and $y\not\in B$. In any way 
            we colour  $A_l$, the colouring rule  extends to a colouring of the closure of $Gx \cup B$ that follows the rule 
             and does not change the colours of $B$. 
             Let $f_l$ be a sequence of such colourings of $Gx \cup B$.  As $Gx$ is a countable set, 
             we can find a   
               point-wise convergent subsequence  of the $f_l$ defining a colouring function $f$ on $Gx \cup B$.  
                We then extend this function $f$ to a pure colouring on the closure of $Gx\cup B$.      \hfill q.e.d.

\section{ Optimality, Stability, and a Bayesian game}

The   colouring rule  $F$ is already formulated as   a problem of local
optimisation according to an objective function.  We can  relax the rules, so that for some given $\delta >0$ it is 
  required   that the colours chosen at all points 
  are within $\delta$ of optimality. We call this {\em point wise $\delta$-optimality}. But that is only one concept 
   of  approximate optimality, that
    at each individual point there is no gain by more than 
    $\delta$ through a different choice of colour.  We seek a slightly broader concept.   
  For each $x\in X$
let $t(x)$ be the possible improvement in the objective function at $x$,
keeping the colouring  for all other $y\neq x$ fixed.  Let $\mu$ be a $G$-invariant
finitely additive extension.  A  colouring  is {\em
  $\gamma$-stable}  if the $\mu$-expectation of $t(x)$ is no more than
$\gamma\geq 0$, meaning that there is no finite disjoint collection
$A_1, \dots , A_n$ of $\mu$ measurable sets such that the objection function
can be improved by at least $t_i$ at all points in~$A_i$
and~$\sum_{i=1} \mu (A_i) t_i $ is greater than $\gamma$.

Another way of
understanding $\gamma$-stability is that~$X$ is a uncountable space of human
society or molecules, and the solution is $\gamma$-stable if the gains from
the individual deviations do not add up to an expectation of $\gamma$.

There are two ways that a measurable colouring  must obey
$\gamma$-stability. First, the set where there is significant divergence from
optimality must be small. Second, where divergence from optimality exists in a
subset of large measure, that divergence must be small (the first concept we presented). That can be formalised
in the following way: if a colouring is $\epsilon \cdot \delta$-stable, then the
subset where it diverges from optimality by more than $\delta$ cannot be of
measure more than $\epsilon$.

 By the continuity of Borel measure, with $p\in \Delta (\{ 1,2,3\} )$ fixed and as $\delta$ goes down to $0$  
    the set of $x$ where   two   rows  of $A_x$ have   expectations within $\delta$ of the optimal choice   falls in Borel measure 
  to $0$. As the  finitely additive measure must be  an extension of the Borel measure, the 
   same is true for the finitely additive measure. 
        This allows the following corollary.    \vskip.2cm

         {\bf Corollary:} For small enough $\gamma$ there is no finitely additive $G$-invariant measure 
          extending the Borel measure
          with a 
  $\gamma$-stable measurable colouring for the colouring rule $F$, and hence there is a 
           paradoxical colouring rule defined by a continuous function.    \vskip.2cm 

{\bf Proof:} Let $\delta>0$ be fixed, and consider the subset $ X_{\delta}$ of $X$ where only one row 
 is   optimal and the other two are not within $\delta$ of being optimal.   
  As $\delta$ 
  goes to zero the probability (according to $m'$) of $X\backslash  X_{\delta}$ goes to zero. 
   Let $\rho(\delta)$ be a function of $\delta$ such  that $\rho (\delta)$  goes to $0$ as $\delta$ goes to $0$  and 
    $\rho (\delta)$ is greater than 
     the probability of $X\backslash  X_{\delta}$ for all positive $\delta$. We break the measurable  colourings
    of $X'$ into two  cases, (1) those involving at least $\frac 1 {100}$ of the space 
      coloured randomly and (2) those with less than $\frac 1{100}$ of the space coloured  randomly. 
      For every  $\epsilon>0$ require that   the
        probability for not  having a pure colour or that  the one  extremal  colour chosen 
         is not optimal is less than $\epsilon$.   We have only 
      to determine  a positive $\epsilon$ and positive $\delta$    small enough  so that 
       $\rho (\delta) < \epsilon$ and $\epsilon$  is small enough to assure, for both cases,  a discrepancy in 
       the measure of the set where purity holds and  $c_1$ is chosen.  For the latter case (2) this is easy,
        $\epsilon  < \frac 1 {1000}$ 
       suffices for a contradiction (as the subset coloured $c_1$ would have to be simultaneously below $\frac 25$ and 
       above $\frac 35$).  
       For the former case (1), due to the continuity of 
        the determinant and therefore also the characteristic polynomial,  following the argument in the proof of Lemma 2 
           the $\frac 1 {100}$ of the space using  random generation  implies 
         the
         existence of a $d>0$ and an $\overline \epsilon>0$ such that   
         the other 
         eigenvalues  (other than $1$) in        the  transition matrix of the colours  
          have norms less than $1-d$ for all  $0<\epsilon< \overline \epsilon $.
           Due to $G$-invariance of any proposed finitely additive  measure,  and assuming $\epsilon < \overline \epsilon$,   
        starting at any distribution on colours and applying $T_1^{-1}$ there is uniform minimal 
          convergence rate to a unique invariant   
          distribution, meaning 
           that at each stage the difference between the distribution and the limit distribution  is no more than 
            $1- d $ times what is was on the previous stage.  As the transitions from the random 
              process determined by $y^e=0$ or $y^e=1$ will remain independent at the various  $y= g_ix$ for independent 
               elements $g_1, \dots , g_k$ and they will  dominate any distributional
                dependency from a small set of measure less 
                than $\epsilon$ of not following the colouring rule $F$,   we can repeat the arguments of 
             Lemma 2, showing that the joint  distributions  of the $g'_i x$ for the seven 
              $g_i'$ are  independent in the limit   as $\epsilon$ goes to $0$.  The same can be done
               for the colour distribution on the whole space, that there is convergence 
                to  independent distributions of colours with $T_1 $ and $T_2$,  using the independence of $T_2 T_1$ and $T_2^2$ 
                 and measure invariance.  
                Due to the constant $d$ and the   existence of  
                  unique invariant distributions, as $\epsilon$ goes to $0$ the convergence   to equalities of  
                 the three  equations with three variables  $y_1, y_2, y_3$  from Lemma 3  implies  that 
                  the unique invariant distributions on colours converges to the 
               same  fixed point  distribution $(y_1, y_2, y_3)= (3 - \sqrt 7, 3-\sqrt 7, 2\sqrt 7 -5)$.  Although the calculations 
               from which one  determines the distribution of colours
                at $x$ from the point $T_1^3 T_2 x= T_2 T_1^3 x$ 
                 don't hold perfectly due to a subset of size  up to $\epsilon$ where the 
                colour rule does not hold, nevertheless 
                   with  $\epsilon$     
           small enough  they show   
           a   persistant $\frac 1 { 5,000}$ discrepancy in 
           the probability given to the colour 
           $c_1$.  
              Having determined positive  
                 $\epsilon< \overline \epsilon $ and $\rho (\delta)< \epsilon$ small enough for 
            both cases, we choose $\delta$ small enough to guarantee  $\rho(\delta)< \epsilon$. The lack of 
             $2 \delta \epsilon$-stability 
             follows. 
             
             Because the matrices $A_x$ change continuously on the  Cantor set $X$ and the two matrices 
              $B_0$ and $B_1$ are defined on disjoint clopen sets (hence together change continuously), and because 
               the $T_1$ and $T_2$ are continuous functions, 
               two points close to each other are close in terms of the consequences of colour choice.  With the colouring rule $F$ defined through the optimisation of a continuous function, 
                optimised at the correspondence $F$, and with uniform continuity of the optimising functions 
                 (from  $X'$ and $C$  compact),
                   any sequence of 
               colouring rule correspondences $F_1, F_2, \dots$   that approximate   the colouring rule  $F$ are also   approximating its 
                its optimisation. With the correspondence $F$ non-empty, upper-semi-continuous and convex valued, for every 
                 sequence of positive $\gamma_1, \gamma _2, \dots$ converging to $0$ 
                there is a sequence of continuous functions  $f_i:X' \times C^2\rightarrow C$ 
                 that approximate  the correspondence $F$ and  the satisfaction of $f_i$ implies 
                 point wise $\gamma_i$ optimality with respect to the correspondence $F$.  
                 By the above,  there is a positive 
                 $\overline \gamma$  where satisfaction of $f_i$ through measurable colouring  is  no longer possible 
                  when  $\gamma_i < \overline \gamma$. Therefore 
                      we get eventually paradoxical colouring
                   rules defined by   continuous 
                 functions.    
             \hfill q.e.d.
          
  We don't use the full force of $\gamma$-stability in showing that there are  paradoxical colouring rules from continuous 
   functions. However there is one applications of $\gamma$-stability that does use that the probability of significant 
    deviation from optimality is limited in probability. 
    
    An important part of economic theory is the study of incomplete information. The idea is that some economic agent has some information that the others do not have, and this private information has to be used carefully to 
     that player's advantage.  It could be only one player with private information or it could be all the players. 
       Often this  situation can  be modelled as a {\em Bayesian}  game.

     The connection to the above colouring rule $F$ is that there is a Bayesian game  played on the same 
      probability space $X'$ for which local 
      optimising behaviour by a player is equivalent to 
      satisfaction of the colouring rule $F$ at an appropriate  point and    
       equilibrium behaviour is equivalent to 
       satisfaction of the  colouring rule $F$ almost everywhere.  Furthermore, a 
         $\gamma$-equilibrium of the game is equivalent to 
        a colouring with the 
        $\gamma$-stability property.

Our interest in paradoxical colouring rules came originally from game theory,
from the desire to show that {\bf all}, not just some, equilibria of a game are
not measurable.  R. Simon \cite{S} showed that there is a Bayesian game which had no
Borel measurable equilibria, though it had non-measurable equilibria. The
infinite dihedral group, an amenable group, acted on the equilibria in a way
that prevented any equilibrium from being measurable. R. Simon and
G. Tomkowicz~\cite{ST} showed that there is a Bayesian game with non-measurable
equilibria but no Borel measurable $\epsilon$-equilibrium for small enough
positive $\epsilon$. That construction involved the action of a non-amenable
semi-group.

A few words are necessary concerning the way a Bayesian game is played. There is a probability space $(\Omega, {\cal F}, m)$; 
 nature chooses a point $x$ in the space $\Omega$  according to the  probability distribution $m$ defined 
  on a sigma algebra ${\cal F}$. There are 
  two approaches  to defining the information, strategies and payoffs of a player $j$. 
  
   In one approach, for each player 
  $j$ there is
   a sigma algebra  ${\cal F}_j$ 
    smaller than ${\cal F}$ such that  the strategy of player $j$ is a function  measurable 
     with respect to ${\cal F}_j$.  
     
      With the other approach  
    we assume  that a  player $j$ has a partition
  ${\cal P}_j$ of 
  the space $\Omega$ -- if nature chooses some $x\in \Omega$   the player $j$  learns that nature's 
  choice lies in the $B\in {\cal P}_j$ such that $x\in B$. If two points $x,y$ belong to the same partition member $B$  then  player $j$ cannot distinguish 
   between $x$ and $y$ and must act identically at $x$ and $y$.

     The different approaches result in different 
       ways to understand what is the strategy and  payoff  of a player. 
       
        With the measurable function 
      approach the  evaluation of a player's  
       strategy  is determined 
       by the strategies of the other players doing the same, and because all their sigma algebras 
        are contained in ${\cal F}$, the evaluation 
         of the result goes through  
          the probability distribution $m$. We call this the Harsanyi approach.
          
           The partition  approach 
          we call the Bayesian approach. With  the Bayesian approach    each player $j$ has 
       a probability distribution on each set in ${\cal P}_j$ and evaluates its actions according to the actions of other 
        players inside of the appropriate member of ${\cal P}_j$. Notice  that the Bayesian approach doesn't really 
         need a probability measure $m$ for the whole space, though we include it because we want to link up 
          the two approaches.   
        
         Though the two approaches are different, they can be related. 
         If  every member of ${\cal P}_j$ 
    is in ${\cal F}$,  
    we  may   move from the Bayesian approach to the Harsanyi approach.     
   A sigma algebra  ${\cal F}_j$ for player $j$ is defined  
   in the canonical way; a set $A\in {\cal F}$  is in ${\cal F}_j$  if and only for every set $B \in {\cal P}_j$
    $A\cap B$ is either $B$ or the empty set.  
          If additionally the player $j$'s probability distributions on each  
         member of ${\cal P}_j$   form a regular conditional 
       probability  with respect to  $m$ and the ${\cal F}_j$ so generated above, we complete the move to the 
        Harsanyi approach.  
     
    Throughout we assume that the payoff of each  player is 
      affine with respect to changes in any one  player's strategy, both with the Harsanyi measurable perspective and with 
     the local Bayesian  perspective.  One can consider more complex payoff structures, but 
       even in the  most trivial information structures  the existence of an equilibrium is not guaranteed 
       when optimality doesn't occur in a convex set.

        Both approaches to what defines a strategy and a payoff 
          have their strengths and weaknesses. The Bayesian approach is more inclusive because 
         it does not require that strategies are measurable. But in general, the Bayesian approach is 
          more problematic. 
         With the Bayesian approach, an evaluation of an action by player $j$ in some 
         $A\in {\cal P}_j$  may be  
         impossible because within $A$ the strategies of the other players may not be measurable with respect to 
          the probability distribution player $j$ has in the set $A$.  But if each member in ${\cal P}_j$ is finite 
           there is not a problem.  Also if  
           there is sufficient structure to the collection ${\cal P}_j$, and we have not yet determined 
            the local probability distributions for each player,  we can determine a probability distribution 
            on each member  $B$ in ${\cal P}_j$ 
           as a  regular conditional   probability  with  
            support on $B$ 
              (see \cite {B}),  and therefore  add a  link  between the   Hansanyi approach 
               to the Bayesian approach.
               
               When the Harsanyi and Bayesian approaches are linked by a regular conditional probability, 
                the importance of measure preserving  invariance to finitely additive extensions can be observed. 
   A  regular conditional probability on some $A\in {\cal P}_j$   must respect (almost everywhere)  any  measure preserving transformations taking 
    place     within the set $A$. 
    By this we mean that if $B$ is a measurable subset   in $A\in {\cal P}_j$ such that $T^{-1} (B)$ is also contained in $A$,   
     then (almost everywhere) the  regular conditional probability  for that player 
     at that set must give $B$ and $T^{-1} (B)$ the same 
     measure.  If a  finitely additive measure on $\Omega$ that extend the original measure $m$ 
      doesn't respect  the local probability distributions of the players,  the game is distorted and  the players' interests 
        would be  no longer represented. 
        If these local probability distributions 
        are defined by measure preserving transformations then it makes sense 
         that the finitely additive measure must keep those transformations 
          measure preserving.

                 An ergodic game, (full definition in \cite{S}), is  
                       one where the most important properties  is that  for each player $j$ each member 
                        of ${\cal P}_j$ is finite and the player's local belief  at each such member 
                         of ${\cal P}_j$ form  a  regular conditional probability.  With  ergodic 
                         games, neither is the  Bayesian approach nor the Harsanyi approaches  problematic.
                            Strictly speaking, the  Bayesian game we present below 
                          is not ergodic because the partition members  of ${\cal P}_j$ for one of the players 
                            are not finite.
                              We  define 
                           a {\em quasi-ergodic} game to have all the same properties of an ergodic game with 
                            the relaxation that  some  members $B$ of  ${\cal P}_j$ may be infinite, 
                             however for every such infinite  $B$ there  a  finite subset $B'\subseteq B$ 
                             such that  inside the set $B\backslash B'$   the player $j$ has no influence over its payoff. Our  Bayesian 
                              game  described below is quasi-ergodic.   
            
            The difference between the two approaches, their different types of strategies and evaluations, gives an added depth to 
             optimisation and stability.  In~\cite {S}
and~\cite{ST} we defined a Harsanyi $\epsilon$-equilibrium for a positive
$\epsilon$: all players in a Harsanyi  $\epsilon$-equilibrium choose measurable
strategies with respect to their sigma algebras as defined above 
 and there is no  measurable deviation by some player to another
measurable strategy resulting in an expected gain of more than $\epsilon$ in  global 
 evaluation. But there is  another 
type of equilibria, the Bayesian. A Bayesian $\epsilon$-equilibrium is a way  for each 
 player to play  that is $\epsilon$-optimal for each set 
 in its partition with respect to its  local    probability distribution on that 
  set. Because the Bayesian equilibrium concept does not require measurable strategies (only that strategies are constant 
    for a player on each set in its partition), there can be Bayesian equilibria where there are no Harsanyi equilibria. 
    This is true for the example in \cite{S}, which is also  an ergodic game.

  There is an added complication to the relation between Harsanyi and Bayesian equilibria when moving 
  to approximate equilibria. 
   Hellman \cite{H} showed that there is  a two person ergodic game without
 a Borel-measurable Bayesian $\epsilon$-equilibria for sufficient small positive $\epsilon$.  For a
positive $\epsilon$,   a
Harsanyi $\epsilon$-equilibrium can be much easier to find than a Harsanyi equilibrium.  
 If the Bayesian game  is defined with an amenable structure (for example through the actions 
 of an amenable group or semi-group)  there will be 
  a Harsanyi $\epsilon$-equilibrium for every $\epsilon >0$ even though there may be no Harsanyi equilibrium \cite {HL} 
   (as  happens with Hellman's example). This is because a Harsanyi
$\epsilon$-equilibrium could employ a very small set 
 where the deviation from local $\epsilon$-equilibrium is
significant, for example of measure
less than $\frac {\epsilon} B$ where the deviation in payoff optimality can be no more than $B$. By  performing 
    this deviation  the measurable behaviour elsewhere of 
    a Harsanyi $\epsilon$-equilibrium could be supported.   
    The existence of such a small set and its role in supporting a Harsanyi
     $\epsilon$-equilibrium  is 
    related to Folner's condition for amenability.

  When moving to finitely additive measures, there are structures to a Bayesian game that enable  paradoxical decompositions.
   To understand this, we take the Bayesian approach and an ergodic game where 
    the information sets of the players are  defined as the orbits of finite groups  that generate a  non-amenable group $G$ and whose elements are measure preserving.  
 There are two intermediate  levels between the whole probability space and individual points in that space. One level 
  is the beliefs of the players as defined by the partitions ${\cal P}_j$ and the probability distributions on each 
   partition member.  As  all sets in all ${\cal P}_j$ are finite,
     this level is very close to the individual points. 
    The other and  higher level  (involving larger sets) 
     is the collection of subsets that the players know in common, 
   the meet partition  of the ${\cal P}_j$, the largest partition smaller than each of the ${\cal P}_j$.    
  The subset that 
        the players know in common, the meet partition, is  the orbit of  $G$. There  may be  
         no    probability distribution  supported on this set, not even a finitely additive one, that is $G$-invariant.
          With the Bayesian approach, the game is played out on  
           these orbits of $G$.   The partitions forming each player's knowledge may be 
            countably generated (the result of countably many refinements of finite partitions of the space) while 
             the meet partition may fail to be  countably generated. This could frustrate any attempt to create Hansanyi
              equilibria 
              from Bayesian  equilibria, including  the broader context of finitely additive measures.  
 Herein  lies the special contribution 
          of Bayesian games to measure theoretic paradoxes.  Nothing  
             pathological about   the information and payoff 
              structures of the individual  players (from the finite group actions) is necessary for the  game to be  paradoxical 
                with respect to equilibria and  finitely additive measures.    
     
The following two  player  Bayesian game has equilibria, yet   fails to have   measurable  
 $\epsilon$-equilibria for sufficiently small $\epsilon>0$, by which we mean 
   measurable  with respect to  any finitely additive extension 
  of the Borel measure that is invariant with respect to the semi-group
   used to define the information structure of the game.

  The most important connection between the information structure of a game  and colouring rules is stated above, 
   that if $x,y$ belong to the same information set of a player, then that player must behave identically at $x$ and $y$. 
    It is this transmission of behaviour over an overlapping system of partitions defined by ergodic operators  
      that connects a  
       colouring rule to equilibrium behaviour.

    We use the same space $X'= X \times   \{ a,b,c\}$ as above. 
    We define two  overlapping  partitions of $\Omega= X'$ corresponding to  two players  $I$ and $II$.   
    
    The information sets of Player $I$ are the sets of the form  $\{ (x,a) \} \cup T_1^{-1} (x) \times \{ b\} \cup 
    T_2^{-1} (x) \times \{ c\}$. The information sets of Player $II$  are the sets of size three 
      of the form $\{ x\} \times \{ a,b,c\}$. 
       Player $I$ considers each of these three sets equally likely. 
        Likewise Player $II$ consider each of the three points equally likely.

          The payoffs to Player $II$  take place   at the   points  $(x,b)$ and $(x,c)$, but with a separate 
           analysis for these two points. Player $II$ has nine actions to choose from, the set of $(j,k)$ corresponding to 
            $j=1,2,3$ and $k=1,2,3$. 
            The payoff for Player $II$ is the sum of the payoffs resulting from  $(x,b)$ and $(x,c)$; 
             one applies the matrix 
             $A_y$ for $y= T_1T_2T_1x$  to the behaviour of Player $I$ at $(x,b)$  
             and the matrix  
             $A_z$  for $z= T_2T_1T_2x$  to the behaviour of Player $I$ at $(x,c)$.  
        The  maximising of  the sum for Player $II$ is accomplished independently. 
         If $p\in \Delta (\{ 1,2,3\})$ is chosen by Player $I$ at $(x,b)$ and  $q\in \Delta (\{ 1,2,3\})$ is chosen 
          by Player $I$ at $(x,c)$ (usually different because they come from different information sets of Player $I$), 
           and Player $II$ chooses the distribution $q_{j,k} \in \Delta (\{ 1,2,3\} \times \{ 1,2,3\})$ 
          then the payoff for Player $II$ is $r A_y p^t + s A_z q^t$, where $y,z$ are defined as above, 
           $r= \sum_k q_{j,k} $ is  Player $II$'s marginal distribution on the first coordinate and 
           $s=\sum_j q_{j,k}$ is Player $II$'s marginal distribution 
            on the second coordinate. 
            
               The payoff to Player $I$ takes place entirely 
             at the point labeled $a$. We would like to say that it follows 
              the three dimensional matrix of the colouring rule $F$, however 
               Player $II$  may not choose an independent distribution.  Nevertheless we can still define a payoff for 
                Player $I$ according to these matrices. We replace the 
                 above evaluation  $\sum_{j=1}^3 \sum _{k=1}^3 p_j q_k b_{i,j,k}$ by 
                  the evaluation $\sum_{j=1}^3 \sum _{k=1}^3 q_{j,k}  b_{i,j,k}$ for the $i$th row, 
                   where $q_{j,k}$ is the probability that Player $II$ chooses the combination  of $j$ with $k$. 
                   
                   The 
                    same argument for purity applies, that if the players are optimising then they are choosing almost everywhere 
                     pure strategies, meaning that Player $I$ puts all weight on only one row and 
                      Player $II$ puts all weight on only  one combination of a $j$ with a $k$.

                         Notice how the causation of players' actions follows through the space. However
Player $I$ behaves at a point $(x,a)$ this behaviour gets translated identically to
$T_1 ^{-1}x\times  \{ b\} $ and to $T_2^{ -1} x\times  \{ c\}$. Player $II$ follows suit 
                  at $T_1 ^{-1}x\times  \{ b\} $ and at $T_2^{ -1} x\times  \{ c\}$
with an attempt to copy Player $I$?s action. However these points belong in uncountably many different 
                   information sets of Player $II$, and the structure of the matrices $A_x$ imply that almost everywhere 
                    the response is pure.  When these various  actions of Player $II$ happen at a 
 common point $(y, a)$ (translated from both $ (y, c) $ and $(y, b)$), Player $I$ responds accordingly
and  transmits this response further  to $T_1^{ -1} (y)\times \{  b\}$ and $T_2^{ -1} (y)\times \{  c\}$.  
The structure of the informations sets mirrors precisely the structure of descendants used to define the colouring rule $F$. When the  behaviour of Player $II$ is pure, the difference in the payoffs of Player $I$ between the game and the 
optimisation of the colouring rule $F$ disappears. With Player $II$'s independent evaluation of the two matrices and the use of marginals to evaluate Player $II$'s payoff, there is no difference to the optimising of  payoffs by Player $II$  and 
 the choosing of optimal rows  that  
 define the colouring rule $F$. Once we have marginalised in
  probability the set where the players are not choosing pure strategies, 
the argument that this game does not have  
 measurable  $\gamma$-equilibria for sufficiently small positive $\gamma$  is the same as that 
                     for the lack of  $\gamma$-stability of any measurable colouring,

     \section {Conclusion} 
     We conjecture  that paradoxical probabilistic colouring 
     rules exist when defined  with group action.         The difficulty seems to lie with the analysis of a stochastic process 
         that isn't  reducible  to a combinatorial argument via purity.   
                                  
                                    \vskip.2cm 
                                  
                                   Does every colouring satisfying our colouring rule $F$, or every  colouring satisfying  a paradoxical 
                                     colouring rule, imply 
                                    the existence of a  measurably $G$-paradoxical decomposition (\cite{ST1}) using  sets of the sigma 
                                      algebra generated by 
   the colour classes, the Borel sets, and the action  of the semi-group generated by  the measure preserving transformations 
    defining the descendants?
     \vskip.2cm 
                
                                                  There is a  problem with applying our above colouring rule $F$ to the group action context. We could     revise the colouring rule to one  on   $X= \{ 0,1\} ^G$ with the group    $G= C_2 * C_3$, 
                            the context of the Hausdorff paradox. 
                            The problem is that the colouring of each point could be so cleverly balanced 
                              as to allow for 
                               the non-purity  of colours throughout the countable orbits of $G$. Putting  these various 
                               orbits together 
                                may still result in a failure of finitely additive 
                                  measurability, but the  argument for this, if true,   is opaque. Surrounding a  three-cycle 
                                with points  coloured purely  
                                 does not stop the colouring of the three-cycle with non-pure colours,
                                  as by Brouwer's Fixed Point Theorem  there 
                                 would  always be 
                                  a  colouring, not necessarily pure,  satisfying the rule inside the three cycle. One 
                                  would have to argue that such a colouring inside the cycle 
                                   must  influence the colouring of the three points 
                                   surrounding that cycle, and as a consequence other  cycles,   such  
                                   that completion of the colouring throughout the orbit, in combination with  other orbits, 
                                      would not be possible in a measurable way. Such an argument 
                                        is plausible, however  seems 
                                       very difficult. 
                  \vskip.2cm 
                  
               The above use of the relation $T_1^3 T_2 = T_2 T_1^3$ was introduced due to the trouble caused 
                by the Brouwer Fixed Point Theorem. Assume there is a probabilistic colouring rule 
                  where the probability space $X$ 
                 is a Cantor set and the  colour is determined 
                 by elements $g_1, \dots, g_k$ of a semi-group of measure preserving  $G$ consisting of 
                  only  non-invertible elements generated by $T_1, \dots , T_m$.   
                   Assume further that 
                  the rule  forces the purity of colours almost 
                 everywhere with the purity of the $g_i x$ forcing 
                  the purity of $x$,   and the rule breaks down into finitely many parts defined 
                  by finitely  many clopen 
                   sets $A_1, \dots , A_l$ that partition $X$.  Let $H$ be the finite subset of $G$ such that 
                    $H=\{ h \ | \ g_i = h j \mbox { for some } j\in G\}$, with $|H| = n$.   
                  We can consider all the possible joint distributions on the colours of 
                  $h y$ for all the $h\in H$ and  $y= T_i x$ for   $i=1, \dots , m $
                   conditioned on the membership 
                   of  the $y$ in the different $ A_j$. If there are relations between the $T_i$, some  joint distribution starting 
                    at some $y=T_ix$ could conflict with another starting at some $y' = T_j x$, and therefore 
                     this approach could lead nowhere.    
                    However given no relation between the $T_i$,
                     the colouring rule  generates in a continuous way (from  purity  and the fixed probabilities 
                      for each of the $A_j$)  a joint distribution on  the colours of 
                    $h x$  for all the $h\in H$  conditioned on $x\in A_j$ for the various 
                     $j=1,\dots, l$. The  Euclidean space 
                      on which we apply Brouwer's fixed point theorem would have dimension  $l \cdot (r^n-1)$ where $r$ 
                       is the number of extremal colours. Though that dimension could be very large, it 
                        is still finite and we 
                         could  apply Brouwer's fixed point 
                     theorem to get  
                      a fixed point of joint distributions, each  conditioned on membership in the $A_1, \dots, A_l$. We 
                      would be robbed of our best tool to demonstrate the lack of an invariant measure. 
                   Is there a way around the Brouwer fixed point theorem that allows for continuous paradoxical colouring 
                    rules while maintaining free generation by non-invertible generations? 
                      \vskip.2cm 
                    
                     Another approach is to single out one non-invertible semi-group element $T$ and  interpret $T^{-1}$ as 
                      the  passage of one unit of time. The idea is that with every passage of a unit of time 
                       there are uncountable variations for continuation, with some conditional 
                        probability distribution governing these variations.
                     In this context one could study how  colouring develops over time. One could use a  paradoxical  colouring 
                      rule for which static  satisfaction   implies finitely additive non-measurability (assuming semi-group invariance), 
                       and a large area satisfying the rule 
                        enjoys some form of relative stability with respect to the passage of time. Is it possible to define a 
                         paradoxical colouring rule this way as something that with high probability, with respect to some 
                           concept of non-measurable  random starts, 
                          will move toward   paradoxical structures of some form? 
                          On the one hand, 
                         local obedience to the rule could be  self perpetuating. On the other hand,   
                         non-amenability of the  structure means  that there is a large  boundary to any rule obedient area 
                          with the potential to destroy that  configuration quickly. 
                           With our above example of a paradoxical colouring rule, 
                            we show that if there is global satisfaction of the rule, it cannot be 
                            invariant measurable. We don't provide an understanding of how that paradoxical structure could 
                             come into existence. 
                \vskip.2cm 
                
              We conjecture that
      there are two-player ergodic  Bayesian games without finitely additive approximate equilibria. 
The discovery of such a Bayesian game 
          would answer all the open problems in the conclusion of \cite {ST}. 
  \vskip.2cm 
                  
   In \cite{ST1} we ask the following: if for a  finite sequence of sets $A_1, A_2, \dots, A_n$  there is no 
     $G$-invariant measure for which all the sets are measurable 
     (where $G$ is the semi-group acting on the probability space) , does
      this imply that  one of these sets  is {\em absolutely} non-measurable, meaning 
                           that there is no finitely additive $G$-invariant measure such that this set is measurable? This is a very relevant 
                            question to game theory. We could easily redefine the above Bayesian game to be a three player game, 
                             with Player $II$ divided into two different players, one evaluating only at points of the form $(x,b)$ 
                              and other one evaluating only at points of the form $(x,c)$.   
         For Player $I$,  who evaluates the 
         colour at $(x,a)$, we could replace the payoff of $1$ with the  the payoff of $2$ whereever that results   
          from the choice of  
          the colour $c_1$. We could do that same for the  player evaluating at $(x,b)$   for the colour 
           $c_2$   and for the other player when the colour is   $c_3$ at an $(x,c)$.    
            In this way, given a positive answer to the above question, the game would become 
            one for which non-measurable 
             equilibria do exist, but none that gives an expected  payoff of any kind  for at least one of the players 
             (which we  consider to be based on some   private  global evaluation not shared  
               with the other players).   
             \vskip.2cm 
             
          One could raise an objection with our above  Bayesian game and our conclusion of no $G$-invariant equilibria.
            Although  with respect to any finitely 
           additive measure that is $G$-invariant there are no measurable 
            strategies in equilibrium, nevertheless a finitely additive measure that  makes all 
            subsets measurable could be the basis for any equilibrium of the game (for which we have proved existence with
             Theorem 2). Indeed Player $I$ receives a non-zero payoff at only one point in each of  its information sets and 
              Player $II$ 
              evaluates actions independently at two different points in each of its corresponding information sets. 
               One could argue that    the local 
              optimisation 
              process remains intact which switching to such a finitely additive measure.  The problem with 
              this perspective is that the   resulting global expected payoff for a player 
              from such a measure  would bear no relation to an integration over the space of  the expected payoffs as 
               perceived locally by the players, 
              since  the conditional probability distributions on their information sets re-confirm the measure 
              invariance of the generators of the semi-group; therefore the probability distributions 
               on their information sets  would have to change when 
               passing to a measure that allowed the global expected payoffs to be measurable.  In the context of game theory, 
               the paradox of our Bayesian 
               game is that one can have equilibria or one can have  payoffs as player 
                expectations of measurable functions, but not both simultaneously. 
                   \vskip.2cm 
                
                There is another answer to the above objection, observed from  
                 a simple change to the payoffs of the game, a change that keeps the game quasi-ergodic.
                   The lack of a measurable equilibrium from a $G$-invariant measure 
                  centres on a discrepancy  
                 of at least $\frac 1 {5,000} $
                 for  what global probability should be  given to  actions corresponding to 
                 the colour $c_1$. For Player $I$ its   information sets 
                  give  one point $(x,a) $  one-third  probability and one-third each  to 
                    uncountable sets of the form $T_1^{-1} x\times  \{ b\} $ and $T_2^{-1} x\times  \{ c\} $, with 
                     the payoff of $0$ for all actions at those two sets. 
                    Let $M$ be a very large positive number. 
                   For the $(x,a)$ point the payoff for  Player $I$ for the action  corresponding to the colour $c_2$ could be increased 
                    for all  combinations of columns   
                    by $M$ and likewise the action  corresponding to  $c_2$ decreased by $M$ uniformly 
                      at the set $T_1^{-1} x\times  \{ b\} $. At the same point, the payoff for Player $I$ for the action corresponding to 
                         $c_3$ could be increased  for all combinations of columns  
                        by $M$ and likewise the action corresponding to $c_3$ 
                        decreased by $M$ uniformly at the set  $T_2^{-1} \times \{ c\}$. 
                      Likewise we do something similar  
                     for Player $II$; at $(x,b)$ we increase the payoff for action corresponding to  $c_2$ for all columns by $2M$ and 
                      decrease it by $M$ uniformly  at the  other two points $(x,a)$ and $(x,c)$ and 
                       do the same for Player $II$ for the action corresponding to  $c_3$  at $(x,c)$ and the other two points.   
                       Notice that nothing is changed strategically as long as the players do not lose their 
                       local probability distributions on their information sets. 
                       Now the definition of an equilibrium of the game must
                        respect the players' local probability evaluations of their information sets, and these in turn 
                         support the measure preserving properties of $T_1$, $T_2$, 
                          and the $S_3$  (that permutes  the $(x,a)$ ,$(x,b)$ and $(x,c)$). With $M$ large enough, any finitely additive 
                           measure for which the strategies of the players are measurable  
                             would distort the  
                       local evaluations of the players sufficiently to violate the equilibrium property. 
                  \vskip.2cm

                 We could consider the following three example, which here 
                  are based on group action. Let $G$ be the group freely generated by $\tau$ and $\sigma$ with 
                   $\sigma^2=e$ and $\tau ^n = e$ for $n\geq 3$. There are $n$ colours $c_1, \dots , c_n$, represented modulo 
                    $n$. 
                     \vskip.2cm 
                     
                   Rule 1: 
                 
                 If $x^e=0$: 
                 
                (1) if  $ \tau^{-1} x$ is coloured with $c_i$, $i\not= n$,   and $\sigma x$ is not coloured $c_1$ then 
                 $x$ is coloured $c_{i}$, 
                 
                 (2) if $\tau^{-1} x$ is coloured with $c_i$ and either $i=n$ or $\sigma x$ is coloured $c_1$ then 
                  $x$ is coloured $c_{i+1}$. 
                  
                  $x^e=1$: if $ \tau ^{-1} x$ is coloured with $c_i$ then $x$ is coloured $c_{i+1}$.

                    \vskip.2cm

                    Rule 2:   
                    
                                   If $x^e=1$
                 and  $ \tau^{-1} x$ is coloured with $c_i$ then colour $x$ with $c_{i+1}$,

                  Given   $x^e=0$ and $\tau^{-1} x $ coloured $c_i$:
                  
                  (1) if  $\sigma x$ is not coloured $c_1$ then colour $x$ with  $c_1$, 
                  
                  (2) if $\sigma x$ is coloured $c_1$  then colour $x$ with  $c_{i+1}$. 
                  \vskip.2cm

                    Rule 3:   
                    
                 The number $n$ is $3$.

                  Given  $\tau^{-1} x $ is coloured $c_i$ and  $\tau x$ is coloured $c_j$:
                  
                 (1) if $i=j$ then colour $x$ with $c_{i+1}$, 
                  
                  (2) if  $i=j+1$  $i+1\not= 1$ and $\sigma x$ is coloured $c_1$ then colour $x$ with $c_{i+1}$, 
                  
                  (3)  if  $i=j+1$  $i+1= 1$ and $\sigma x$ is not coloured $c_1$ then colour $x$ with  $c_{i+1}$, 
                  
                  (4) in all other cases, colour $x$ with $c_j$.

                  \vskip.2cm 

                   All three  rules are very interesting, but we are not sure where they could lead in terms 
                    of establishing probabilistic  paradoxical colouring rules or other interesting colouring structures. It is relatively 
                     easy to show that Rule 1 is paradoxical for $n\geq 5$ (from a discrepancy for the probability of the colour $c_1$),
                      but it is not clear if it could be the basis of a probabilistic 
                      paradoxical rule.  
                     The third rule, defined in \cite {MB} and proven there to be paradoxical, would need to be revised to 
                      a rule dependent on location in order for it to be the basis of a probabilistic colouring rule.


\begin{thebibliography}{HD}

\bibitem{B} D. Blackwell and L Dubins, \emph {Existence and Non-Existence of Proper, Regular, Conditional Distributions },
 The Annals of Probability, Vol. 3, No. 5, (1975), pp. 741-752.

\bibitem{MB} M. Bounds, \emph {M.Sc. thesis}, London School of Economics, 2020. 

\bibitem{H} Z. Hellman, \emph {A Game with no Bayesian Approximate Equilibria},
 Journal of Economic Theory, 153(C) (2014), pp. 138-151.

\bibitem{HL} Z. Hellman and  J. Yehuda , \emph {Equilibria Existence in Bayesian Games: Climbing the Countable Borel Equivalence Relation Hierarchy}, Working Papers 2020- 15, Business School - Economics, University of Glasgow.



\bibitem{S} R. S. Simon, Games of Incomplete Information, Ergodic 
Theory, and the Measurability of Equilibria,
 {\em Israel J. Math.}, 138, 1,  (2003) pp. 73-92. 

 \bibitem{ST} R. S. Simon and G. Tomkowicz,
      \emph {A Bayesian Game without $\epsilon$-equilibria},
       Israel J. of Math. 227 (2018), 215-231.   
   

\bibitem{ST1} R. S. Simon and G. Tomkowicz,
      \emph {Paradoxical decompositions and finitary colouring rules},
       Preprint.  

\end{thebibliography}
\end{document}